\documentclass{sig-alternate}

\usepackage{color}
\usepackage{pdflscape}
\usepackage[table]{xcolor}
\usepackage{algorithm}
\usepackage{algorithmicx}
\usepackage{algpseudocode}
\usepackage{subfig}

\definecolor{lightgray}{RGB}{217,224,250}

\newdef{definition}{Definition}
\newdef{example}{Example}

\graphicspath{{figures/}}

\usepackage[
  pass,
]{geometry}

\begin{document}

\conferenceinfo{GECCO'14} {July 12-16, 2014, Vancouver, BC, Canada.}
\CopyrightYear{2014}
\crdata{TBA}
\clubpenalty=10000
\widowpenalty = 10000

\title{Certified Global Minima for a Benchmark of Difficult Optimization
Problems}

\numberofauthors{4}
\author{
~ 
~ \alignauthor
~ Charlie Vanaret \\
~ \affaddr{ENAC\titlenote{Ecole Nationale de l'Aviation Civile}, IRIT\titlenote{Institut de Recherche en Informatique de Toulouse}}\\
~ \affaddr{7 Avenue Edouard Belin}\\
~ \affaddr{31055 Toulouse Cedex 4}\\
~ \affaddr{France}\\
~ \email{vanaret@cena.fr}
~ 
~ \alignauthor
~ Jean-Baptiste Gotteland \\
~ \affaddr{ENAC, IRIT}\\
~ \affaddr{7 Avenue Edouard Belin}\\
~ \affaddr{31055 Toulouse Cedex 4}\\
~ \affaddr{France}\\
~ \email{gottelan@cena.fr}
~ 
~ \alignauthor
~ Nicolas Durand \\
~ \affaddr{ENAC, IRIT}\\
~ \affaddr{7 Avenue Edouard Belin}\\
~ \affaddr{31055 Toulouse Cedex 4}\\
~ \affaddr{France}\\
~ \email{durand@cena.fr}
~ 
~ \and
~ \alignauthor
~ Jean-Marc Alliot \\
~ \affaddr{IRIT}\\
~ \affaddr{118 Route de Narbonne}\\
~ \affaddr{31062 Toulouse Cedex 9}\\
~ \affaddr{France}\\
~ \email{jean-marc.alliot@irit.fr}
}

\maketitle
\begin{abstract}
We provide the global optimization community with new optimality proofs
for six deceptive benchmark functions (five bound-constrained functions
and one nonlinearly constrained problem). These highly multimodal
nonlinear test problems are among the most challenging benchmark
functions for global optimization solvers; some have not been solved
even with approximate methods.

The global optima that we report have been numerically certified using
Charibde (Vanaret \textit{et al.}, 2013), a hybrid algorithm that
combines an evolutionary algorithm and interval-based methods. While
metaheuristics generally solve large problems and provide sufficiently
good solutions with limited computation capacity, exact methods are
deemed unsuitable for difficult multimodal optimization problems. The
achievement of new optimality results by Charibde demonstrates that
reconciling stochastic algorithms and numerical analysis methods is a
step forward into handling problems that were up to now considered
unsolvable.

We also provide a comparison with state-of-the-art solvers based on
mathematical programming methods and population-based metaheuristics,
and show that Charibde, in addition to being reliable, is highly
competitive with the best solvers on the given test functions.
\end{abstract}




\section{Introduction}
Numerical solvers usually embed advanced methods to tackle nonlinear
optimization problems. Stochastic methods, in particular evolutionary
algorithms, handle large problems and provide sufficiently good
solutions with limited computation capacity, but may easily get trapped
in local minima. On the other hand, local and global (exhaustive)
deterministic methods may guarantee local or global optimality, but are
often limited by the size or the nonlinearity of the problems and may
suffer from numerical approximations.

In \cite{Vanaret2013}, a new reliable hybrid solver named Charibde has
been introduced to reconcile stochastic methods and numerical analysis
methods. An evolutionary algorithm and an interval-based algorithm are
combined in a cooperative framework: the two methods run in parallel and
cooperate by exchanging the best known upper bound of the global minimum
and the best current solution. The contribution of this paper is the
achievement of new certified optimality results by Charibde for six
highly multimodal nonlinear test functions, for which no or few results
were available. We also compare Charibde with state-of-the art solvers
including mathematical programming methods, population-based
metaheuristics and spatial branch and bound. Charibde proves to be
highly competitive with the best solvers on the given test functions,
while being fully reliable.

Charibde is presented in Section \ref{sec:charibde}. It is evaluated on
a benchmark of difficult test functions given in Section
\ref{sec:functions}. State-of-the-art solvers to which Charibde is
compared are described in Section \ref{sec:solvers}. Numerical results,
including proofs of optimality, values of global minima and
corresponding solutions, are provided and discussed in Section
\ref{sec:results}.

\section{Charibde: a Rigorous Solver}\label{sec:charibde}
The rigorous nonlinear solver Charibde was introduced by Vanaret
\textit{et al.} \cite{Vanaret2013}, building on an original idea by
Alliot \textit{et al.} \cite{Alliot2012}. It combines the efficiency of
a Differential Evolution algorithm and the reliability of interval
computations to discard more efficiently subspaces of the search-space
that cannot contain a global minimizer.
We introduce interval computations and interval-based methods in
Sections \ref{sec:ia} and \ref{sec:ibc}. Details on the implementation
of the Differential Evolution algorithm are given in Section
\ref{sec:de}. Finally, the cooperation scheme of Charibde is explained
in Section \ref{sec:cooperation}.

\subsection{Interval Analysis}\label{sec:ia}
\textbf{Interval Analysis} (IA) is a method of numerical analysis
introduced by Moore \cite{Moore1966} to bound rounding errors in
floating-point computations. Real numbers that are not representable on
a computer are enclosed within intervals with floating-point bounds.
Each numerical computation is safely carried out by using outward
rounding.

\begin{definition}
An interval $X = [\underline{X}, \overline{X}]$ is the set
$\{x \in \mathbb{R} ~|~ \underline{X} \le x \le \overline{X} \}$.
We denote by $m(X) = \frac{1}{2}(\underline{X} + \overline{X})$ its
midpoint. $\mathbb{IR}$ is the set of all intervals.
A box $\mathbf{X} = (X_1, \ldots, X_n)$ is an interval vector.
We note $m(\mathbf{X}) = (m(X_1), \ldots, m(X_n))$ its midpoint.
In the following, capital letters represent interval quantities
(interval $X$) and bold letters represent vectors (box $\mathbf{X}$,
vector $\mathbf{x}$).
\end{definition}

\textbf{Interval arithmetic} defines the interval counterparts of
real-valued operators ($\{+, -, \times, /\}$) and elementary functions
($\exp$, $\cos$, $\ldots$). For example, $[a, b] + [c, d] = [a
+_{\downarrow} c, b +_{\uparrow} d]$ and $\exp([a, b]) =
[\exp_{\downarrow}(a), \exp_{\uparrow}(b)]$, where $\cdot_{\downarrow}$
(resp. $\cdot_{\uparrow}$) denotes downward (resp. upward) rounding.

\begin{definition}
Let $f: \mathbb{R}^n \rightarrow \mathbb{R}$ be a real-valued function.
$F : \mathbb{IR}^n \rightarrow \mathbb{IR}$ is an \textit{interval
extension} of $f$ if
\begin{equation*}
\begin{aligned}[]
\forall \mathbf{X} \in \mathbb{IR}^n,  f(\mathbf{X}) =
\{f(\mathbf{x}) ~|~ \mathbf{x} \in \mathbf{X} \} \subset F(\mathbf{X}) \\
\forall (\mathbf{X}, \mathbf{Y}) \in \mathbb{IR}^n,
\mathbf{X} \subset \mathbf{Y} \Rightarrow F(\mathbf{X}) \subset F(\mathbf{Y}) &
\end{aligned}
\end{equation*}
The \textit{natural interval extension} $F_N$ is obtained by replacing
elementary operations in $f$ with their interval extensions.
\end{definition}

\textbf{Dependency} is the main source of overestimation when using
interval computations: multiple occurrences of a same variable are
considered as different variables. For example, the interval evaluation
of $f(x) = x^2 - 2x$ over the interval $[1, 4]$ yields $F_N([1, 4]) =
[-7, 14]$, which crudely overestimates the exact range $f([1, 4]) = [-1,
8]$. However, an appropriate rewriting of the syntactic expression of
$f$ may reduce or overcome dependency: if $f$ is continuous inside a box,
$F_N$ yields the optimal range when each variable occurs only once in
its expression. Completing the square in the expression of $f$ provides
the optimal syntactic expression $g(x) = (x - 1)^2 - 1$. Then $G_N([1,
4]) = [-1, 8] = f([1, 4])$.

\subsection{Interval-based Techniques}\label{sec:ibc}

\subsubsection{Interval Branch and Bound Algorithms}
Interval Branch and Bound algorithms (IB\&B) exploit the
conservative properties of interval extensions to rigorously bound
global optima of numerical optimization problems~\cite{Hansen1992}.
The method consists in splitting the initial search-space into subspaces
(branching) on which an interval extension is evaluated (bounding). By
keeping track of the best upper bound $\tilde{f}$ of the global minimum
$f^*$, boxes that certainly do not contain a global minimizer are
discarded (example \ref{ex:ibb}). Remaining boxes are stored to be
processed at a later stage until the desired precision $\varepsilon$ is
reached. The process is repeated until all boxes have been processed.
Convergence certifies that $\tilde{f} -
f^* < \varepsilon$, even in the presence of rounding errors. However, the
exponential complexity of IB\&B hinders the speed of convergence on
large problems.

\begin{example}
\label{ex:ibb}
Consider the problem $\min\limits_{x \in X} f(x) = x^4 - 4x^2$ over
$X = [-1, 4]$. Then
$F_N([-1, 4]) = [-64, 256] \supset [-4, 192] \\= f([-1, 4])$. The
floating-point evaluation $f(1) = -3$ provides an upper bound $\tilde{f}$
of $f^*$. Evaluating $F_N$ on the subinterval $[3, 4]$ reduces the
overestimation induced by dependency: $F_N([3, 4]) = [17, 220] \supset
[45, 192] = f([3, 4])$. Because $\forall x \in [3, 4], f(x) \ge 17 >
\tilde{f} = -3 \ge f^*$, the interval $[3, 4]$ cannot contain a global
minimizer and can be safely discarded.
\end{example}

\subsubsection{Interval Contraction}
Propagating the (in)equality constraints of the problem, as well as the
constraints $f \le \tilde{f}$ and $\nabla f = 0$, may narrow the
domains of the variables or prove that a subdomain of the search-space
cannot contain a global minimizer.

Stemming from the IA and Interval Constraint Programming communities,
filtering/contraction algorithms \cite{ChabertJaulin2009} narrow the
bounds of the variables without loss of solutions. Standard contraction
algorithms generally integrate a filtering procedure into a fixed-point
algorithm. \texttt{HC4} \cite{Benhamou1999} handles one constraint after
the other and performs the optimal contraction w.r.t. to a constraint if
variables occur only once in its expression. \texttt{Box}
\cite{VanHentenryck1997} narrows one variable after the other w.r.t. all
constraints, using an interval version of Newton's method. \texttt{Mohc}
\cite{Araya2010} exploits the monotonicity of the constraints to enhance
contraction of \texttt{HC4} and interval Newton.

The interval-based algorithm embedded in Charibde follows an Interval
Branch and Contract (IB\&C) scheme (algorithm \ref{alg:ibc}) that
interleaves steps of bisection and filtering. We note $\mathcal{L}$ the
priority queue in which the remaining boxes are stored, $\varepsilon$
the desired precision and $\tilde{x}$ the best known solution, such that
$\overline{F(\tilde{x})} = \tilde{f}$.

\begin{algorithm}[h!]
\caption{Interval Branch and Contract framework}
\label{alg:ibc}
\begin{algorithmic}
\State $\tilde{f} \leftarrow +\infty$
\Comment best found upper bound
\State $\mathcal{L} \leftarrow \{ \mathbf{X}_0 \}$
\Comment priority queue of boxes to process
\Repeat
	\State Extract a box $\mathbf{X}$ from $\mathcal{L}$
	\Comment selection rule
	\State Compute $F(\mathbf{X})$
	\Comment bounding rule
	\If{$\mathbf{X}$ cannot be eliminated}
	\Comment cut-off test
		\State Contract$(\mathbf{X}, \tilde{f})$
		\Comment{filtering algorithms}
		\State Compute $F(m(\mathbf{X}))$ to update $\tilde{f}$
		\Comment{midpoint test}
		\State Bisect $\mathbf{X}$ into $\mathbf{X}_1$ and $\mathbf{X}_2$
		\Comment branching rule
		\State Store $\mathbf{X}_1$ and $\mathbf{X}_2$ in $\mathcal{L}$
	\EndIf
\Until{$\mathcal{L} = \varnothing$}
\State \Return $(\tilde{f}, \tilde{x})$
\end{algorithmic}
\end{algorithm}


\subsection{Differential Evolution}\label{sec:de}
Differential Evolution (DE) is an Evolutionary Algorithm that combines
the coordinates of existing individuals with a particular probability to
generate new potential solutions \cite{StornPrice1997}. It was embedded
within Charibde for its ability to solve extremely difficult
optimization problems, while having few control parameters.

We denote by $NP$ the population size, $W > 0$ the weighting factor and
$CR \in [0, 1]$ the crossover rate. For each individual $\mathbf{x}$ of
the population, three other individuals $\mathbf{u}$, $\mathbf{v}$ and
$\mathbf{w}$, all different and different from $\mathbf{x}$, are
randomly picked in the population. The newly generated individual
$\mathbf{y} = (y_1, \ldots, y_j, \ldots, y_n)$ is computed as follows:
\begin{equation}
	y_j =
	\begin{cases}
	u_j + W \times (v_j - w_j) & \text{ if } j = R \text{ or } r_j < CR \\
	x_j & \text{ otherwise}
	\end{cases}
\end{equation}
where $R$ is a random index in $\{1, \ldots, n\}$ ensuring that at least one
component of $\mathbf{y}$ differs from that of $\mathbf{x}$, and
$r_j$ is a random number uniformly distributed in $[0, 1]$, picked for
each $y_j$. $\mathbf{y}$ replaces $\mathbf{x}$ in the population if
$f(\mathbf{y}) < f(\mathbf{x})$.

The following advanced rules have been implemented in Charibde:
\begin{description}
\item[Boundary constraints:] When a component $y_j$ lies outside the
bounds $[\underline{X_j}, \overline{X_j}]$ of the search-space, the
\textit{bounce-back method} \cite{PriceStornLampinen2006} replaces
$y_j$ with a component that lies between $u_j$ and the admissible bound:
\begin{equation}
\small
	y_j =
	\begin{cases}
	u_j + rand(0, 1) (\overline{X_j} - u_j), & \text{if } y_j > \overline{X_j} \\
	u_j + rand(0, 1) (\underline{X_j} - u_j), & \text{if } y_j < \underline{X_j}
	\end{cases}
\end{equation}

\item[Evaluation:]
Given inequality constraints $\{g_i ~|~ i = 1, \ldots, m\}$, the
evaluation of an individual $\mathbf{x}$ is computed as a triplet
($f_\mathbf{x}$, $n_\mathbf{x}$, $s_\mathbf{x}$), where $f_\mathbf{x}$
is the objective value of $\mathbf{x}$, $n_\mathbf{x}$ the number of violated
constraints and $s_\mathbf{x} = \sum_{i = 1}^m \max(g_i(\mathbf{x}), 0)$.
If at least one of the constraints is violated, the objective value is
not computed

\item[Selection:]
Given the evaluation triplets ($f_\mathbf{x}, n_\mathbf{x},
s_\mathbf{x}$) and ($f_\mathbf{y}, n_\mathbf{y}, s_\mathbf{y}$) of two
candidate solutions $\mathbf{x}$ and $\mathbf{y}$, the best individual
to be kept for the next generation is computed as follows:
\begin{itemize}
\item if $n_\mathbf{x} < n_\mathbf{y}$ or ($n_\mathbf{x} = n_\mathbf{y}
> 0$ and $s_\mathbf{x} < s_\mathbf{y}$) or ($n_\mathbf{x} = n_\mathbf{y}
= 0$ and $f_\mathbf{x} < f_\mathbf{y}$) then $\mathbf{x}$ is kept
\item otherwise, $\mathbf{y}$ replaces $\mathbf{x}$
\end{itemize}
\end{description}

\subsection{Charibde: a Cooperative Algorithm}\label{sec:cooperation}
\begin{figure}[h!]
\centering
\def\svgwidth{\columnwidth}
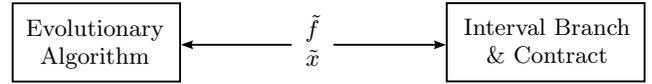
\caption{Basic scheme of cooperation}
\label{fig:cooperation}
\end{figure}

Charibde combines an Interval Branch and Contract algorithm and a
Differential Evolution algorithm in a \textit{cooperative} way: neither
of the algorithms is embedded within the other, but they run in parallel
and exchange bounds and solutions using an MPI implementation (Figure
\ref{fig:cooperation}).

The cooperation scheme boils down to 3 main steps:
\begin{itemize}
\item Whenever the best known DE evaluation is improved, the
best individual $\mathbf{x}_{b}$ is evaluated using IA. The upper bound
of the image $\overline{F(\mathbf{x}_{b})}$ -- an upper bound of the
global minimum -- is sent to the IB\&C thread
\item In the IB\&C algorithm, $\overline{F(\mathbf{x}_{b})}$ is compared
to the current best upper bound $\tilde{f}$. An improvement of the
latter leads to a more efficiently pruning of the subspaces that cannot
contain a (feasible) global minimizer
\item Whenever the evaluation of the center $m(\mathbf{X})$ of a
box improves $\tilde{f}$, the individual $m(\mathbf{X})$ replaces the
worst individual of DE, thus preventing premature convergence
\end{itemize}

\section{Benchmark of Test Functions}\label{sec:functions}
\begin{figure*}[t]
	\centering
	\subfloat[Michalewicz function]{\epsfig{file=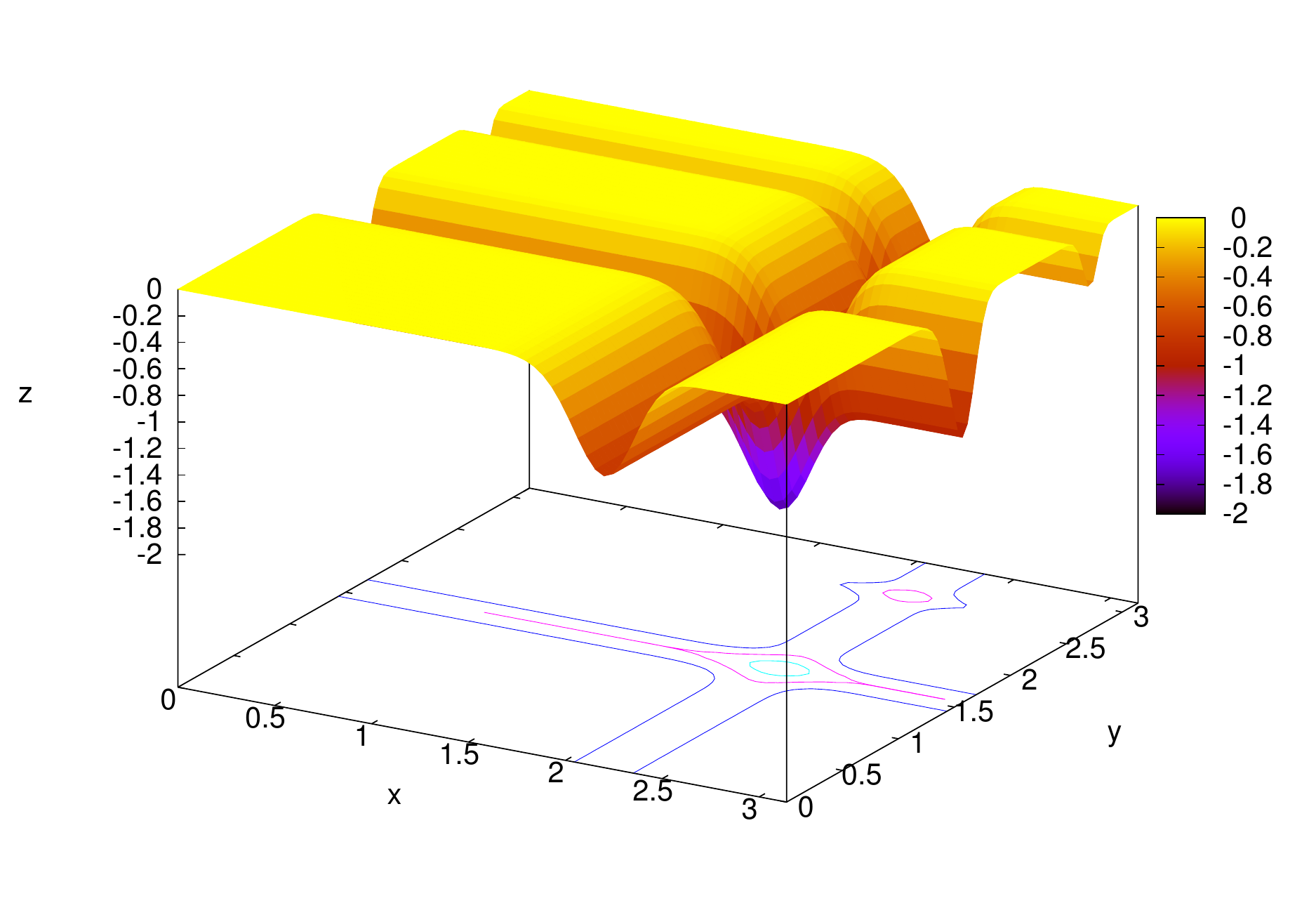, width=5.9cm}}
	\subfloat[Sine Envelope function]{\epsfig{file=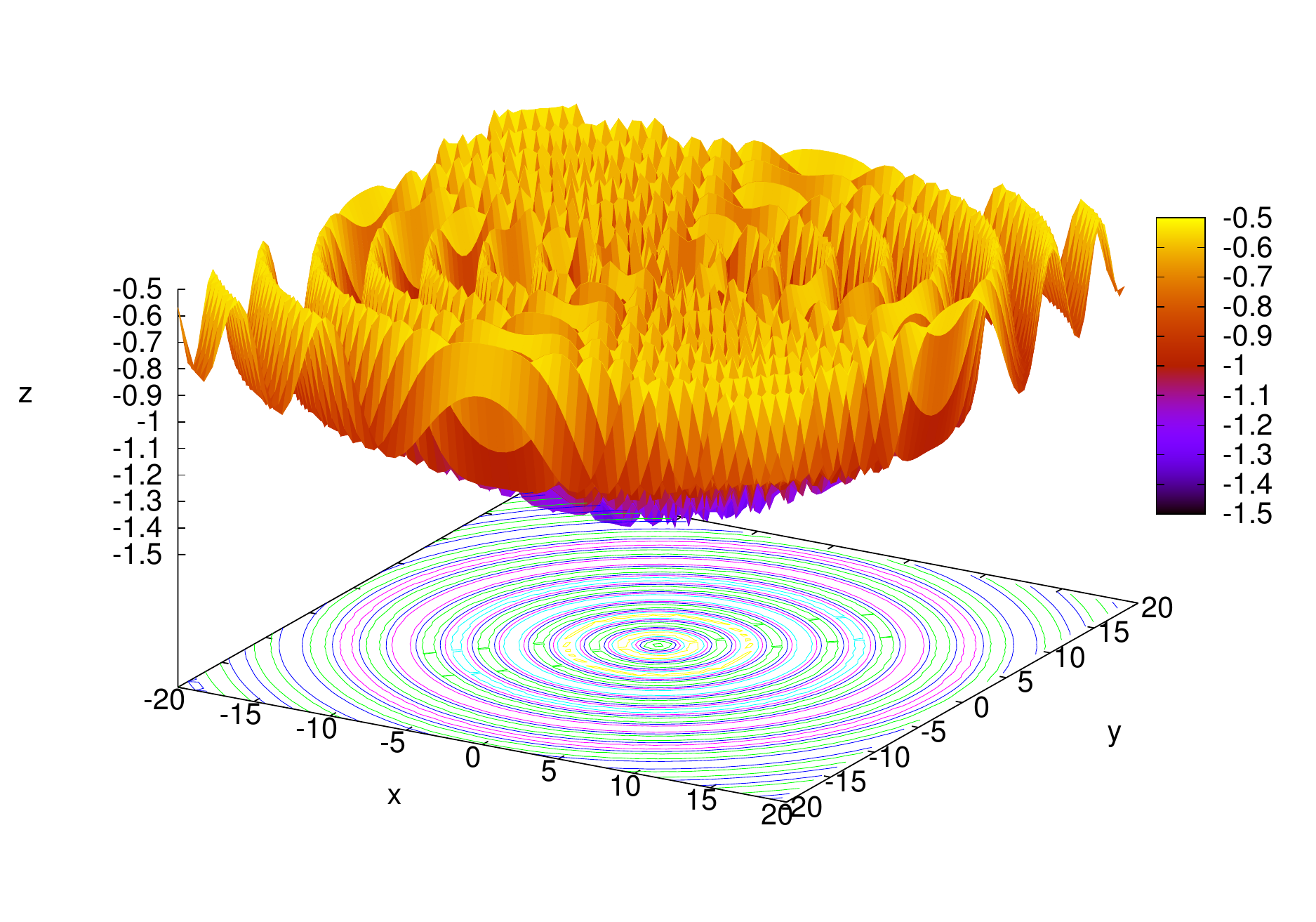, width=5.9cm}}
	\subfloat[Shekel's Foxholes function]{\epsfig{file=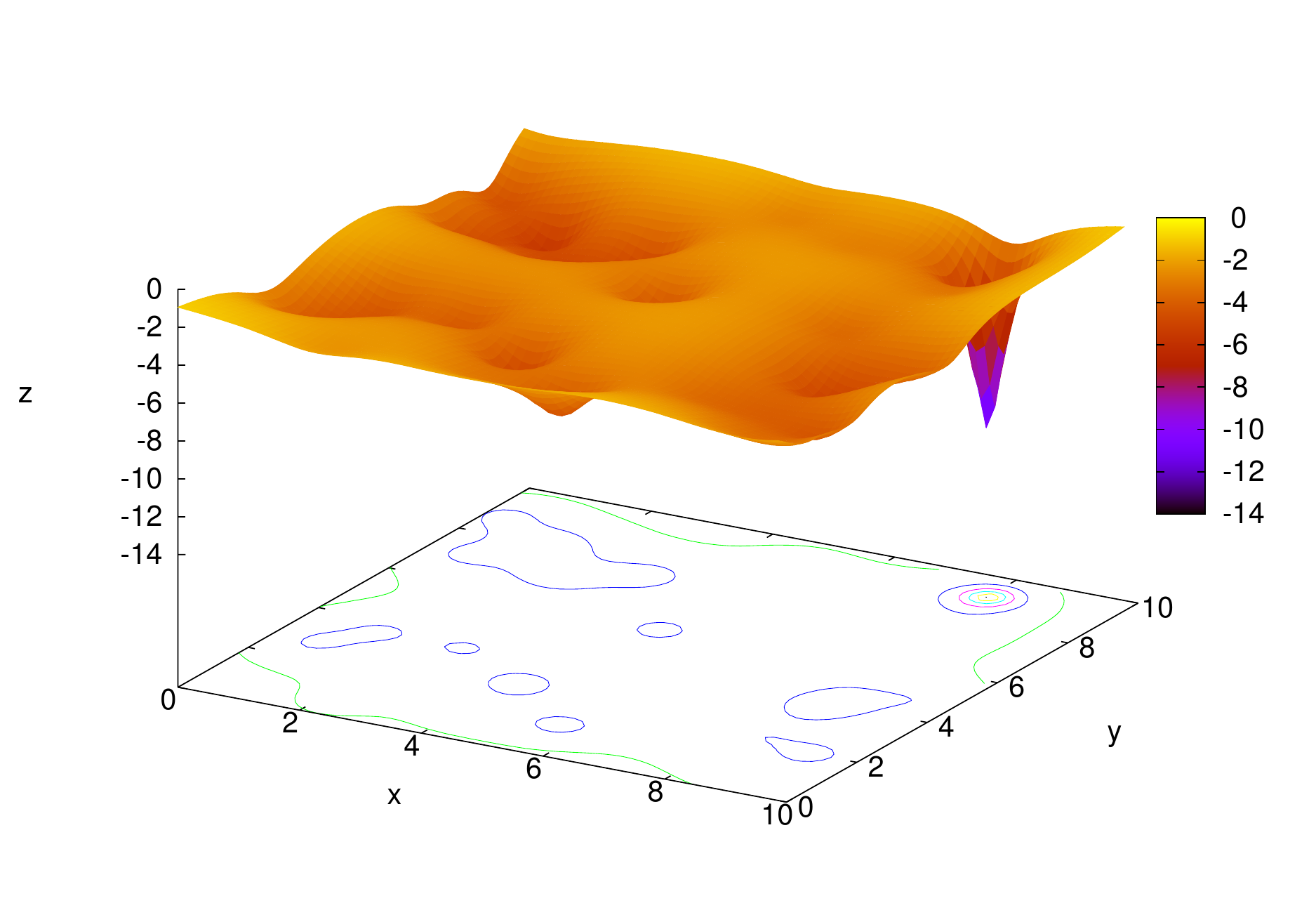, width=5.9cm}} \quad
	\subfloat[Egg Holder function]{\epsfig{file=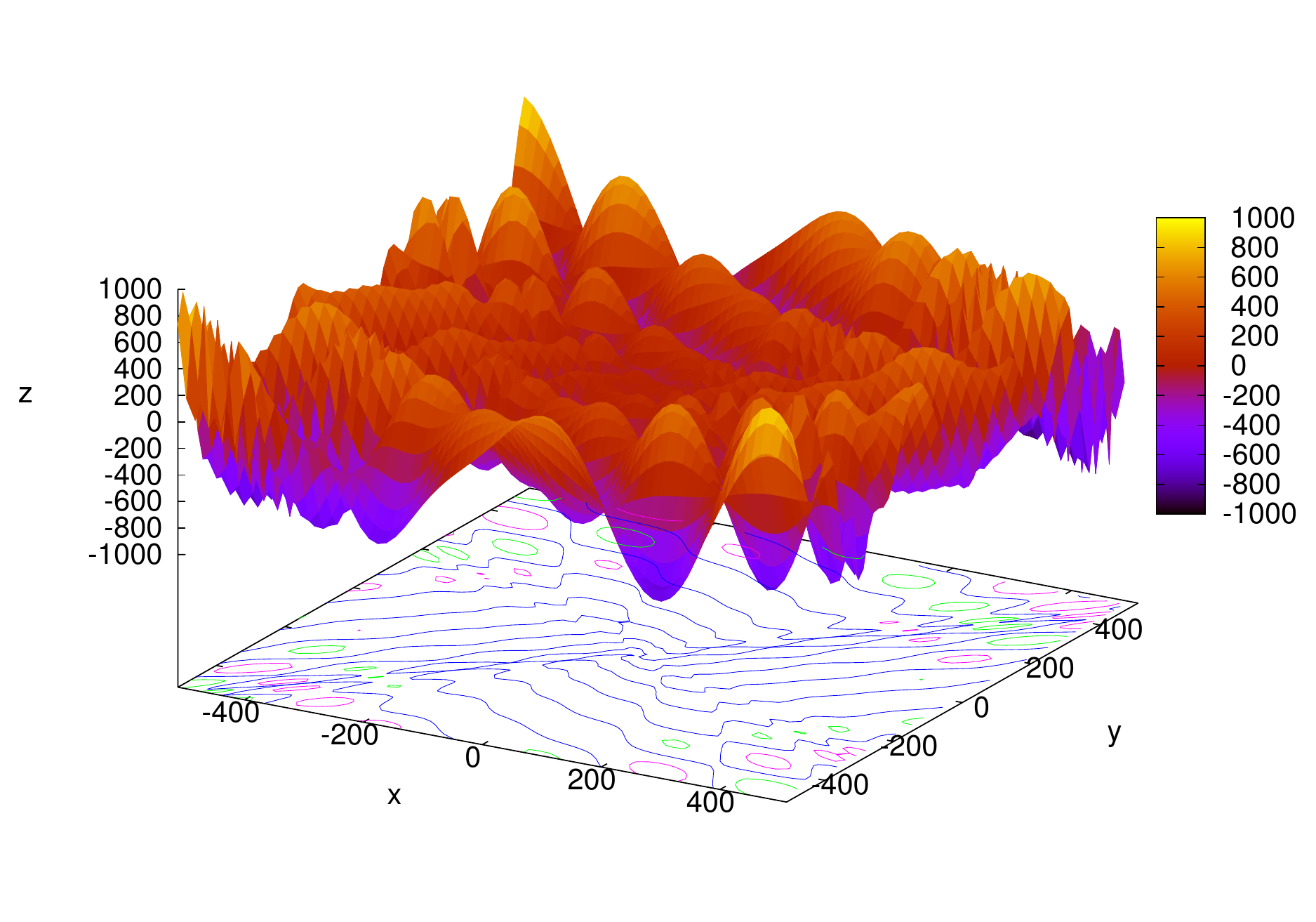, width=5.9cm}}
	\subfloat[Rana's function]{\epsfig{file=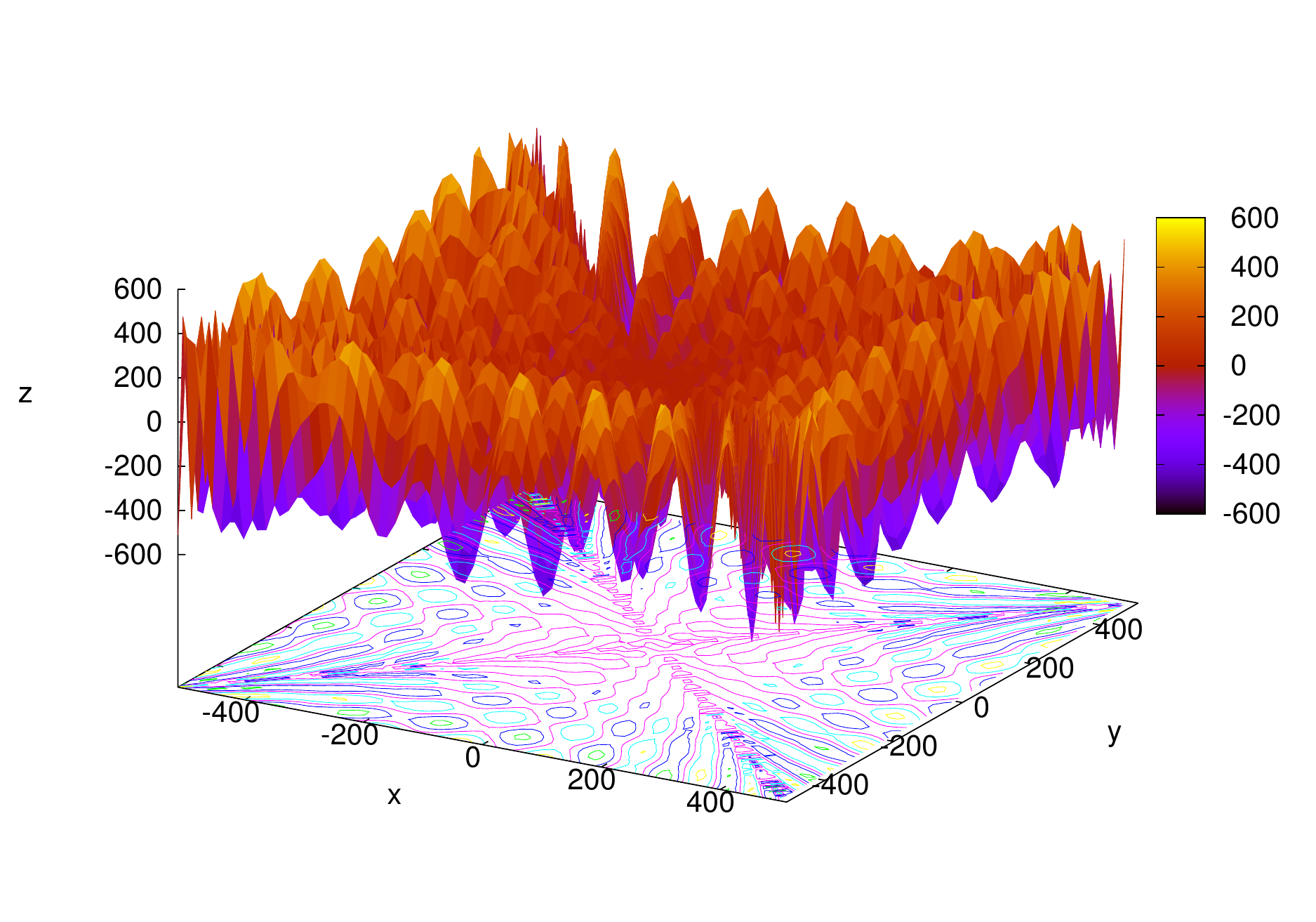, width=5.9cm}}
	\subfloat[Keane's function]{\epsfig{file=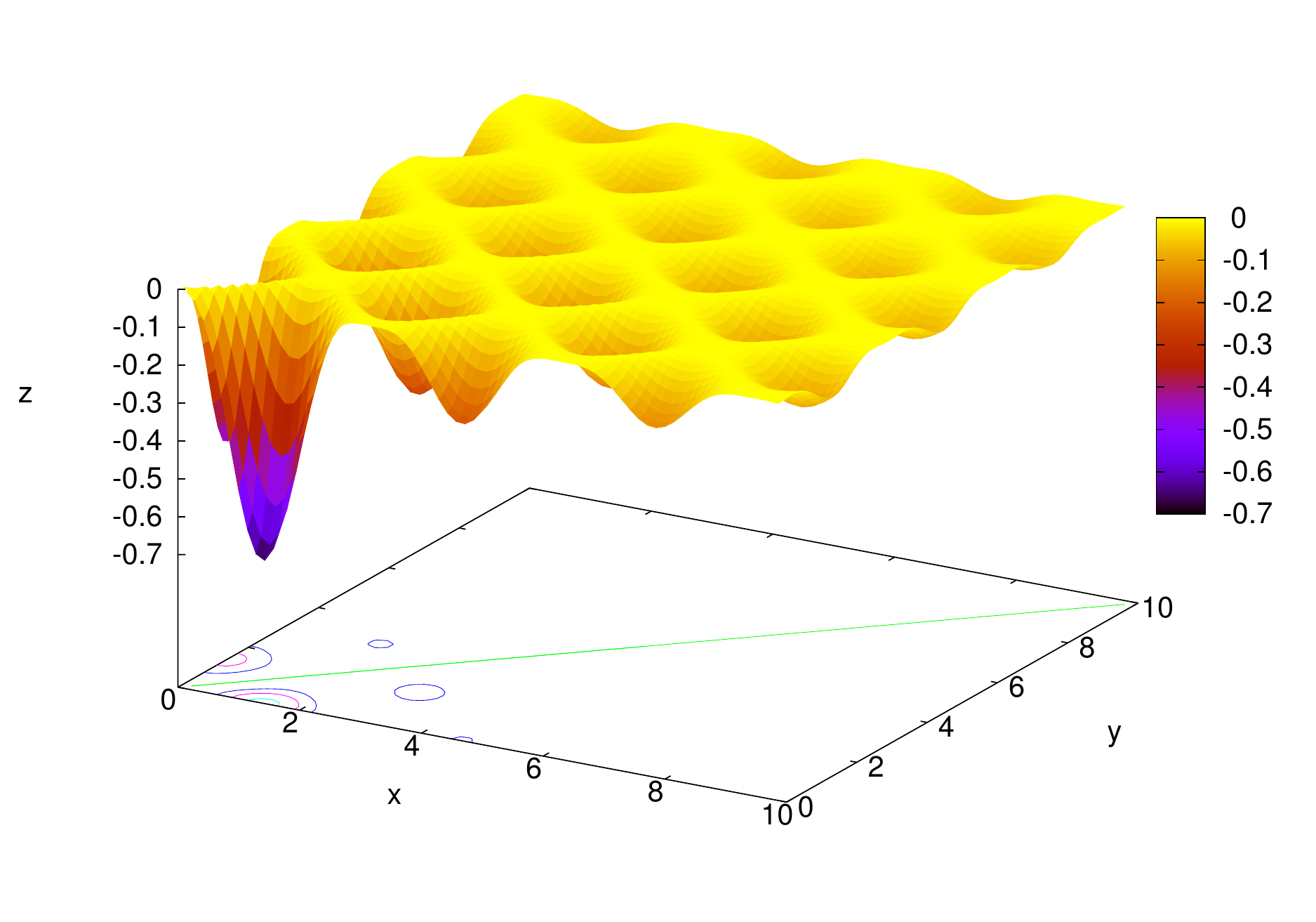, width=5.9cm}}
	\caption{Test functions ($n = 2$)}
	\label{fig:test-functions}
\end{figure*}

The highly multimodal nonlinear test functions considered in this study
can be found in Table \ref{tab:test-functions}. Contrary to standard
test functions that have a global minimum 0 at $(0, \ldots, 0)$
(Griewank function) or have a global minimizer with $n$ identical
components (Schwefel function), we have selected six functions with
nontrivial global minima: \textbf{Michalewicz} function, \textbf{Sine
Envelope Sine Wave} function (shortened to Sine Envelope),
\textbf{Shekel's Foxholes} function \cite{PriceStornLampinen2006},
\textbf{Egg Holder} function \cite{Whitley1996}, \textbf{Rana}'s
function \cite{Whitley1996} and \textbf{Keane}'s function
\cite{Keane1994}. Except for the Michalewicz function, all are
nonseparable. Their surfaces and contour lines can be observed on
Figure \ref{fig:test-functions} for $n = 2$. It illustrates the numerous
local minima and the ruggedness of the functions.

The first inequality constraint of Keane's function describes a
hyperbola in two dimensions and is active at the global minimizer, which
hinders the efficiency of solvers. The second inequality constraint is
linear and is not active at the global minimizer.
The Egg Holder (resp. Rana) function is strongly subject to dependency:
$x_1$ and $x_n$ occur three (resp. five) times in its expression, and
$(x_2, \ldots, x_{n-1})$ occur six (resp. ten) times. Its natural
interval extension therefore produces a large overestimation of the
actual range.

The last three functions (Egg Holder, Rana and Keane) contain absolute
values. $|\cdot|$ is differentiable everywhere except for $x = 0$,
however its subderivative -- generalizing the derivative to
non-differentiable functions -- at $x = 0$ can be computed. Charibde
handles an interval extension proposed by Kearfott \cite{Kearfott1996},
based on the values of its subderivative:
\begin{equation}
	|\cdot|'(X) =
	\begin{cases}
	[-1, -1] & \text{ if } \overline{X} < 0 \\
	[1, 1] & \text{ if } \underline{X} > 0 \\
	[-1, 1] & \text{otherwise}
	\end{cases}
\end{equation}
This expression is used to compute an enclosure of the partial
derivatives of $f$ through automatic differentiation.

\begin{table*}[t]
	\centering
	\scriptsize
	\rowcolors{1}{}{lightgray}
	\caption{Test functions}
	\begin{tabular}{|llc|}
	\hline
	Function & Expression & Domain \\
	\hline
	Michalewicz 	& $-\sum_{i=1}^n \sin(x_i) \left[ \sin(\frac{i x_i^2}{\pi}) \right]^{20}$ & $[0, \pi]^n$ \\
	Sine Envelope 	& $-\sum_{i=1}^{n-1}\left( 0.5 + \frac{\sin^2(\sqrt{x_{i+1}^2+x_{i}^2}-0.5)}{(0.001(x_{i+1}^2+x_{i}^2)+1)^2} \right)$ & $[-100, 100]^n$ \\
	Shekel			& $-\sum_{i=1}^{30} \frac{1}{c_i + \sum_{j=1}^n (x_j - a_{ij})^2}$ & $[0, 10]^n$ \\
	Egg Holder 		& $-\sum_{i=1}^{n-1} \left[ (x_{i+1}+47) \sin \left(\sqrt{|x_{i+1} + 47 + \frac{x_i}{2}|} \right) + x_i \sin \left(\sqrt{|x_i-(x_{i+1}+47)|} \right) \right]$ & $[-512, 512]^n$ \\
	Rana 			& $\sum_{i=1}^{n-1} ( x_i \cos \sqrt{|x_{i+1}+x_i+1|} \sin \sqrt{|x_{i+1}-x_i+1|} + (1+x_{i+1}) \sin \sqrt{|x_{i+1}+x_i+1|} \cos \sqrt{|x_{i+1}-x_i+1|})$ & $[-512, 512]^n$ \\
	Keane 			& $-\frac{|\sum_{i=1}^n \cos^4 x_i - 2\prod_{i=1}^n \cos^2 x_i|}{\sqrt{\sum_{i=1}^n ix_i^2}}$ s.t. $0.75 \le \prod_{i=1}^n x_i$ and $\sum_{i=1}^n x_i \le 7.5n$ & $[0, 10]^n$ \\
	\hline
	\end{tabular}

	\label{tab:test-functions}
\end{table*}

\section{State-of-the-art Solvers}\label{sec:solvers}
The certified global minima obtained by Charibde are compared with
state-of-the-art optimization softwares available on the NEOS
(Network-Enabled Optimization System) server\footnote{
\texttt{http://www.neos-server.org/neos/}}, a free web service for
solving optimization problems in AMPL/GAMS formats. These solvers
include local methods, population-based metaheuristics and spatial
branch and bound: \textbf{Ipopt} implements a primal-dual interior point
algorithm, which uses a filter line search method. \textbf{LOQO} is
based on an infeasible primal-dual interior-point method. \textbf{MINOS}
employs a reduced-gradient method. \textbf{PGAPack} is a Parallel
Genetic Algorithm library. \textbf{PSwarm} combines pattern search and
particle swarm. \textbf{Couenne} uses a reformulation-based
branch-and-bound algorithm for a globally optimum solution. It is a
complete solver, i.e. it performs an exhaustive exploration of the
search-space.

Note: the solver \textbf{BARON} (Branch and Reduce Optimization
Navigator), combining constraint propagation, interval analysis and
duality, is known as one of the most efficient, although unreliable (see
Section \ref{sec:couenne}), solvers for solving nonconvex optimization
problems to global optimality. However, it is not designed to handle
trigonometric functions, and therefore will not be considered in this
paper.

\vspace{2cm}

\section{Numerical Experiments}\label{sec:results}

\subsection{Global Minima Certified by Charibde}
The global minima to the six test functions are given in Table
\ref{tab:minima}. All solutions are given with a precision $\varepsilon
= 10^{-6}$. Should the functions have several global minima, only one is
provided. A reference given next to a global minimum indicates that
a result is available in the literature.

Analyzing the results has brought out identities for finding putative
minima for four of the six test functions. They can be found in Table
\ref{tab:putative-minima}.

\begin{table}[h!]
	\centering
	\rowcolors{1}{}{lightgray}
	\caption{Putative minima}
	\begin{tabular}{|llc|}
	\hline
	Function & Putative minimum & $R^2$ \\
	\hline
	Michalewicz 	& $-0.99864n + 0.30271$ & 0.9999981 \\
	Sine Envelope 	& $-1.49150n + 1.49150$ & 1 \\
	Egg Holder		& $-915.61991n + 862.10466$ & 0.9999950 \\
	Rana's			& $-511.70430n + 511.68714$ & 1 \\
	\hline
	\end{tabular}
	
	\label{tab:putative-minima}
\end{table}

The global minima of the Michalewicz function for up to 75 variables can
be found in Table \ref{tab:micha-minima}. For the sake of conciseness,
the corresponding solutions for up to only 10 variables are given in
Table \ref{tab:minima}.

\begin{table}[h!]
	\centering
	\rowcolors{1}{}{lightgray}
	\caption{Certified global minima of Michalewicz function}
	\begin{tabular}{|cl|cl|}
	\hline
	$n$ & Global minimum & $n$ & Global minimum \\
	\hline
	10 & -9.6601517 \cite{Mishra2006} 	& 45 & -44.6256251 \\
	15 & -14.6464002 					& 50 & -49.6248323 \cite{Mishra2006} \\
	20 & -19.6370136 \cite{Mishra2006} 	& 55 & -54.6240533 \\
	25 & -24.6331947 					& 60 & -59.6231462 \\
	30 & -29.6308839 \cite{Mishra2006} 	& 65 & -64.6226167 \\
	35 & -34.6288550 					& 70 & -69.6222202 \\
	40 & -39.6267489 					& 75 & -74.6218112 \\
	\hline
	\end{tabular}
	
	\label{tab:micha-minima}
\end{table}

To illustrate the key role played by the syntactic expression of the
function when computing with intervals, we have tested two different --
but equivalent -- syntaxes of Rana's function in Charibde. The first
syntax is given in Table \ref{tab:test-functions}. The second syntax is
obtained using the trigonometric identity: $\cos x \sin y = \frac{1}{2}(
\sin(x+y) - \sin(x-y))$. Their impact on the sharpness of the inclusions,
therefore on the convergence of Charibde, can be observed in Table
\ref{tab:rana-rewriting}. The hyphen indicates a computing time greater
than one hour.

\begin{table}[h!]
	\centering
	\rowcolors{1}{}{lightgray}
	\caption{Rana's function: CPU times (s) of convergence in Charibde
	($NP = 70$, $W = 0.7$, $CR = 0.5$)}
	\begin{tabular}{|ccc|}
	\hline
	$n$ & \multicolumn{2}{c|}{CPU time (s)} \\
		& First syntax & Second syntax \\
	\hline
	2 & 0.25 & 0.009 \\
	3 & 6.5 & 0.12 \\
	4 & 254 & 1.45 \\
	5 & - & 18.5 \\
	6 & - & 244 \\
	7 & - & 3300 \\
	\hline
	\end{tabular}
	
	\label{tab:rana-rewriting}
\end{table}

\subsection{Comparison of Solvers}
The comparison of the seven solvers on a particular instance of each
test function can be found in Table \ref{tab:test}. When available, the
number of evaluations of the objective function is given under the
found minimum, otherwise we have mentioned the computing time. The
hyphens in the last column indicate that PGAPack and PSwarm cannot
handle inequality constraints.

Local methods based on mathematical programming (Ipopt, LOQO, MINOS)
usually require few iterations to reach a local minimum from a starting
point. The quality of the minimum generally depends on both the starting
point and the size of the basins of attraction of the function. These
solvers turn out to perform poorly on the considered multimodal problems.

Among the population-based metaheuristics, PGAPack performs consistently
better than PSwarm, yet at a higher cost. We have kept the best result
ouf of five runs, since it is tedious to run an online solver several
times. The default NEOS control parameters have been used; both
algorithms would certainly perform better with suitable control
parameters. The results of Couenne are discussed in Section
\ref{sec:couenne}.

Charibde achieves convergence in finite time on the six problems, with a
numerical certification of optimality. It benefits from the start of
convergence of the DE algorithm that computes a good initial value for
$\tilde{f}$. This allows the IB\&C algorithm to prune more efficiently
subspaces of the search-space. The number of evaluations of the natural
interval extension $F$ has the form $NE_{DE} + NE_{IB\&C}$, where
$NE_{DE}$ is the number of evaluations in the DE algorithm (whenever the
best known solution is improved) and $NE_{IB\&C}$ is the number of
evaluations in the IB\&C algorithm. Note that after converging toward
the global minimizer, the DE thread keeps running as long as the
certification of optimality has not been obtained.

\begin{table*}[t!]
	\centering
	\rowcolors{1}{}{lightgray}
	\caption{State-of-the-art solvers against Charibde on test functions}
	\begin{tabular}{|lcccccc|}
	\hline
	 		& Michalewicz 	& Sine Envelope & Shekel 		& Egg Holder  	& Rana  		& Keane  \\
			& ($n = 50$)	& ($n = 5$)		& ($n = 5$)		& ($n = 5$)		&($n = 5$)		& ($n = 4$) \\
	\hline
	Ipopt 	& -19.773742 	& -5.8351843 	& -1.8296708 	& -3586.3131827 & -75.512076	& -0.2010427 \\
		 	& (167 eval) 	& (24 eval) 	& (15 eval) 	& (8 eval) 		& (16 eval) 	& (8 eval) \\
	LOQO 	& -0.0048572 	& -5.8351843 	& -1.9421143 	& 25.5609238 	& -69.5206		& -0.0983083 \\
		 	& (88 eval) 	& (17 eval) 	& (19 eval) 	& (5924 eval)	& (138 eval)	& (50 eval) \\
	MINOS 	& 0 			& -5.87878 		& -2.5589419 	& -3586.313183 	& -233.592		& -0.2347459 \\
		 	& (3 eval) 		& (38 eval) 	& (3 eval) 		& (3 eval)		& (1 eval)		& (3 eval) \\
	PGAPack & -37.60465 	& -5.569544 	& -1.829452 	& -3010.073 	& -2091.068		& - \\
			& (9582 eval) 	& (9615 eval) 	& (9602 eval) 	& (9626 eval) 	& (9622 eval)	& - \\
	PSwarm 	& -24.38158 	& -5.835182 	& -1.610072 	& -2840.799 	& -1595.056		& - \\
		 	& (2035 eval) 	& (2049 eval) 	& (2008 eval) 	& (2040 eval) 	& (2046 eval)	& - \\
	Couenne & -49.6\underline{19042}		& -5.96\underline{60007} 	& \textbf{-10.4039521} 	& -3719.72\underline{87498} & \textbf{-2046.8320657} & -0.6222\underline{999} \\
			& (265s) 		& (0.4s) 		& (8s) 					& (20.6s) & (20.3s) & (2s) \\
	\hline
	\textbf{Charibde}	& \textbf{-49.6248323} 	& \textbf{-5.9659811}	& \textbf{-10.4039521} 	& \textbf{-3719.7248363} 	& \textbf{-2046.8320657} 	& \textbf{-0.6222810} \\
	CPU time			& 8.4s	 				& 219s					& 0.04s 				& 0.8s	 					& 17.8s 					& 14860s \\
	$f$ evaluations		& 717100				& 227911400				& 8150					& 742600					& 19186150 					& 6351693297 \\
	$F$ evaluations 	& 763 + 409769			& 93 + 21744667			& 28 + 561				& 106 + 82751				& 53 + 1383960 				& 74 + 2036988566 \\
	($NP, W, CR$)		& (50, 0.7, 0)			& (50, 0.7, 0.9) 		& (50, 0.7, 0.9)		& (50, 0.7, 0.4) 			& (50, 0.7, 0.5) 			& (70, 0.7, 0.9) \\
	\hline
	\end{tabular}
	
	\label{tab:test}
\end{table*}

\subsection{Reliability vs Efficiency}\label{sec:couenne}
\begin{figure*}[h!]
	\centering
	\includegraphics[width=1\textwidth]{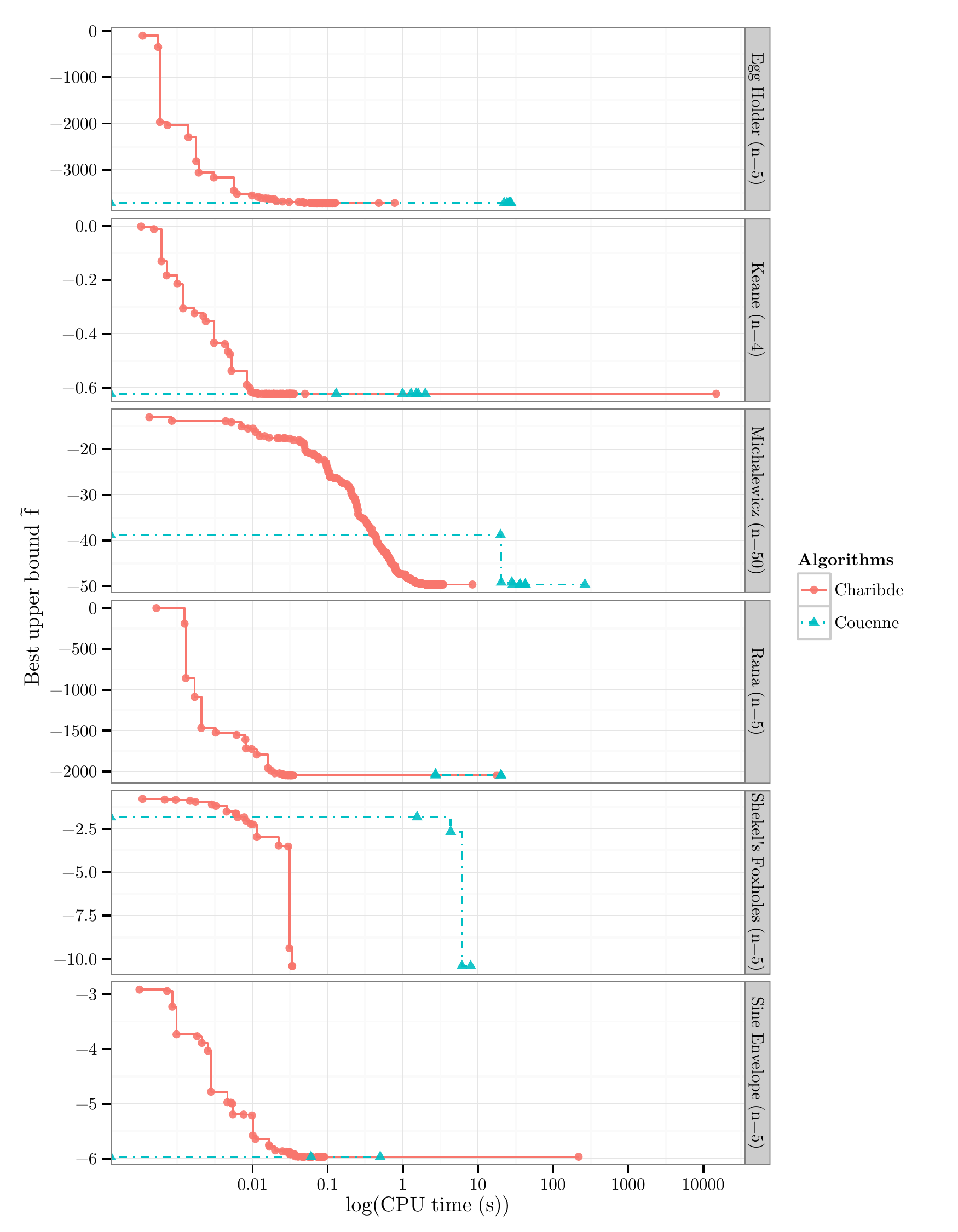}
	\caption{Comparison of Couenne and Charibde (logarithmic x scale)}
	\label{fig:charibde-couenne}
\end{figure*}

On Figure \ref{fig:charibde-couenne}, the evolution of the best known
solution in Charibde is compared with that of Couenne for a particular
instance of each test function. Intermediate times for other solvers
were not available. Note that the scale on the x-axis is logarithmic.
These diagrams show that Charibde is highly competitive against Couenne:
Charibde achieves convergence faster than Couenne on Michalewicz (ratio
31), Shekel (ratio 200), Egg Holder (ratio 25) and Rana's function
(ratio 1.1). Couenne is faster than Charibde on Sine Envelope (ratio
547) and Keane's function (ratio 7430).

Is is however crucial to note that Couenne, while being a complete
solver (the whole search-space is exhaustively processed), is
\textit{not reliable}. This stems from the fact that the under- and
overestimators obtained by linearizing the function may suffer from
numerical approximations. Contrary to interval arithmetic that bounds
rounding errors, mere real-valued linearizations are not conservative
and cannot guarantee the correctness of the result. This problem can
easily be observed in Table \ref{tab:test}: the global minima obtained
by Couenne on Michalewicz, Sine Envelope, Egg Holder and Keane's
function are not correct compared to the certified minima provided by
Charibde. The wrong decimal places are underlined.

\section{Conclusion}
We provided a comparison between Charibde, a cooperative solver that
combines an EA and interval-based methods, and state-of-the-art solvers
(stemming from mathematical programming and population-based
metaheuristics). They were evaluated on a benchmark of nonlinear
multimodal optimization problems among the most challenging:
Michalewicz, Sine Envelope, Shekel's Foxholes, Egg Holder, Rana and
Keane. Charibde proved to be highly competitive with the best solvers,
including Couenne, a complete but unreliable solver based on spatial
branch and bound and linearizations.
We provided new certified global minima for the considered test
functions, as well as the corresponding solutions. They may be used from
now on as references to test stochastic or deterministic optimization
methods.


\bibliographystyle{abbrv}
\bibliography{Vanaret}

\begin{landscape}
\begin{table}
	\rowcolors{1}{}{lightgray}
	\caption{Certified global minima}
	\begin{tabular}{|l|c|c|l|}
	\hline
	Function & $n$ & Global minimum & Corresponding solution \\
	\hline
	Michalewicz	& 2 	& -1.8013034 & (2.202906, 1.570796) \\
				& 3		& -2.7603947 & (2.202906, 1.570796, 1.284992) \\
				& 4		& -3.6988571 & (2.202906, 1.570796, 1.284992, 1.923058) \\
				& 5		& -4.6876582 & (2.202906, 1.570796, 1.284992, 1.923058, 1.720470) \\
				& 6		& -5.6876582 & (2.202906, 1.570796, 1.284992, 1.923058, 1.720470, 1.570796) \\
				& 7		& -6.6808853 & (2.202906, 1.570796, 1.284992, 1.923058, 1.720470, 1.570796, 1.454414) \\
				& 8		& -7.6637574 & (2.202906, 1.570796, 1.284992, 1.923058, 1.720470, 1.570796, 1.454414, 1.756087) \\
				& 9		& -8.6601517 & (2.202906, 1.570796, 1.284992, 1.923058, 1.720470, 1.570796, 1.454414, 1.756087, 1.655717) \\
				& 10	& -9.6601517 & (2.202906, 1.570796, 1.284992, 1.923058, 1.720470, 1.570796, 1.454414, 1.756087, 1.655717, 1.570796) \\

	Sine Envelope 	& 2 & -1.4914953 \cite{Pohl2010} & (-0.086537, 2.064868) \\
					& 3 & -2.9829906 & (1.845281, -0.930648, 1.845281) \\
					& 4 & -4.4744859 & (2.066680, 0.001365, 2.066680, 0.001422) \\
					& 5 & -5.9659811 & (-1.906893, -0.796823, 1.906893, 0.796823, -1.906893) \\
					& 6 & -7.4574764 & (-1.517016, -1.403507, 1.517016, -1.403507, -1.517015, 1.403507) \\

	Shekel	& 2 & -12.1190084 \cite{Koullias2013} & (8.024065, 9.146534) \\
			& 3 & -11.0307623 & (8.024161, 9.150962, 5.113211) \\
			& 4 & -10.4649942 & (8.024876, 9.151655, 5.113888, 7.620843) \\
			& 5 & -10.4039521 \cite{Feng2013} & (8.024917, 9.151728, 5.113927, 7.620861, 4.564085) \\
			& 6 & -10.3621514 & (8.024916, 9.151795, 5.113951, 7.620875, 4.564063, 4.710999) \\
			& 7 & -10.3131505 & (8.024945, 9.151838, 5.113968, 7.620912, 4.564052, 4.711004, 2.996069) \\
			& 8 & -10.2793068 & (8.024947, 9.151874, 5.113979, 7.620933, 4.564043, 4.711003, 2.996055, 6.125980) \\
			& 9 & -10.2288309 & (8.024961, 9.151914, 5.113990, 7.620956, 4.564025, 4.710999, 2.996038, 6.125996, 0.734065) \\
			& 10 & -10.2078768 \cite{Pedroso2007} & (8.024968, 9.151929, 5.113991, 7.620959, 4.564020, 4.711005, 2.996030, 6.125993, 0.734057, 4.981999) \\

	Egg Holder 	& 2 	& -959.6406627 \cite{Oplatkova2008}	& (512, 404.231805) \\
				& 3		& -1888.3213909 & (481.462894, 436.929541, 451.769713) \\
				& 4 	& -2808.1847922 & (482.427433, 432.953312, 446.959624, 460.488762) \\
				& 5 	& -3719.7248363 & (485.589834, 436.123707, 451.083199, 466.431218, 421.958519) \\
				& 6 	& -4625.1447737 & (480.343729, 430.864212, 444.246857, 456.599885, 470.538525, 426.043891) \\
				& 7 	& -5548.9775483 & (483.116792, 438.587598, 453.927920, 470.278609, 425.874994, 441.797326, 455.987180) \\
				& 8 	& -6467.0193267 & (481.138627, 431.661180, 445.281208, 458.080834, 472.765498, 428.316909, 443.566304, 457.526007) \\
				& 9 	& -7376.2797668 & (482.785353, 438.255330, 453.495379, 469.651208, 425.235102, 440.658933, 454.142063, 468.699867, 424.215061) \\
				& 10 	& -8291.2400675 & (480.852413, 431.374221, 444.908694, 457.547223, 471.962527, 427.497291, 442.091345, 455.119420, 469.429312, 424.940608) \\
	Rana 	& 2 & -511.7328819 \cite{Tao2007} & (-488.632577, 512) \\
			& 3 & -1023.4166105 	& (-512, -512, -511.995602) \\
			& 4 & -1535.1243381 	& (-512, -512, -512, -511.995602) \\
			& 5 & -2046.8320657 	& (-512, -512, -512, -512, -511.995602) \\
			& 6 & -2558.5397934 	& (-512, -512, -512, -512, -512, -511.995602) \\
			& 7 & -3070.2475210 	& (-512, -512, -512, -512, -512, -512, -511.995602) \\
	Keane 	& 2 & -0.3649797 \cite{Kang2002} & (1.600860, 0.468498) \\
		 	& 3 & -0.5157855 \cite{Kang2002} & (3.042963, 1.482875, 0.166211) \\
		 	& 4 & -0.6222810 \cite{Kang2002} & (3.065318, 1.531047, 0.405617, 0.393987) \\
	\hline
	\end{tabular}
	\label{tab:minima}
\end{table}
\end{landscape}


\end{document}